\documentclass[11pt,openbib]{article}
\usepackage{amssymb,amsmath,amsfonts,theorem}
\newcommand{\R}{\mathbb{R}}

\newcommand{\gl}{\mathop{\mathrm{\! \, gl}}\nolimits}
\newcommand{\Gl}{\mathop{\mathrm{\! \, Gl}}\nolimits}
\newcommand{\oo}{\mathop{\mathrm{\! \, o}}\nolimits}
\newcommand{\OO}{\mathop{\mathrm{\! \, O}}\nolimits}

\newcommand{\comp}{\raisebox{0pt}{$\scriptstyle\circ \, $}}
\newcommand{\setrule}{\, \rule[-4pt]{.5pt}{13pt}\, }

\newcommand{\spann}{\mathop{\rm span}\nolimits}

\allowdisplaybreaks
\begin{document}

\begin{center}
{\Large \bf Adjoint orbits of the Lie algebra \\
\rule{0pt}{18pt} of the generalized real orthogonal group} \\
\mbox{} \vspace{.05in} \\
Richard Cushman
\end{center}
\footnotetext{printed: \today}
\bigskip 

\section{Adjoint orbits of the generalized real orthogonal group}  

In this section we define the generalized real orthogonal group and 
introduce some concepts needed to classify its adjoint orbits on 
its Lie algebra. \medskip 

First we give a description of the generalized real 
orthogonal group in terms 
of matrices. Let $\widetilde{G}$ be the Gram matrix of a nondegenerate inner product 
$\widetilde{\gamma }$ on ${\R }^n$. 
The \emph{real orthogonal group} 
$\OO ({\R}^n, \widetilde{G})$ on $({\R}^n, \widetilde{G})$ is the set 
of $n\times n$ real matrices $\widetilde{A}$ such that 
${\widetilde{A}}^T\widetilde{G}\widetilde{A} =\widetilde{G}$. It is a Lie 
group whose Lie algebra $\oo ({\R }^n,\widetilde{G})$ is 
$\{ \widetilde{\xi} \in \gl (n, \R) \setrule \, 
{\widetilde{\xi}}^{\, T}\widetilde{G} + \widetilde{G} \widetilde{\xi } =0\} $. Let $\{ e_0, \, e_1, \ldots , e_n \}$ be 
the standard basis for ${\R }^{n+1}$. 
Suppose that $\gamma $ is symmetric bilinear form on ${\R }^{n+1}$ 
whose Gram matrix is $G =${\tiny $ \left( \begin{array}{c|c} 0 & 0 \\ \hline 
\rule{0pt}{7pt} 0 & \widetilde{G} \end{array} \right) $}. The bilinear form 
$\gamma $ is \emph{degenerate} since $\ker G = \spann \{ e_0 \} $. The 
\emph{generalized real orthogonal group} ${\OO ({\R }^{n+1}, G)}_{e_0}$ on 
$({\R }^{n+1}, G)$ is $\{ A \in \OO ({\R }^{n+1}, G) \, \setrule \, A e_0 = e_0 \}$, that is, 
\begin{displaymath}
\left\{ \mbox{\footnotesize $ \left( \begin{array}{c|c} 1 & b^T \\ \hline 0 & 
\rule{0pt}{12pt}\widetilde{A} \end{array} \right) $} \right.  \left. \in \Gl (n+1, \R ) \setrule 
\, b \in {\R }^n, \, \, \widetilde{A} \in \OO ({\R }^n, \widetilde{G})
\right\} . 
\end{displaymath}
${\OO ({\R}^{n+1}, G)}_{e_0}$ is a Lie group with Lie algebra ${\oo ({\R }^{n+1}, G)}_{e_0}$ given by 
\begin{displaymath}
\left\{ \mbox{\footnotesize $\left( \begin{array}{c|c} 0 & b^T \\ \hline
\rule{0pt}{12pt} 0 & \widetilde{\xi } \end{array} \right) $} \right.  \left. \in \gl (n+1, \R ) \setrule \, b\in {\R }^n, \, \, \widetilde{\xi } \in \right. \left. \oo ({\R }^n, \widetilde{G} ) \right\} . 
\end{displaymath} 

We now give a more abstract description of the generalized real orthogonal 
group. Let $\widetilde{V}$ be a finite dimensional real vector space with a 
symmetric bilinear form $\widetilde{\gamma }$, which is  nondegenerate.   
Let $\OO (\widetilde{V}, \widetilde{\gamma})$ be the set of real linear mappings $\widetilde{A}$ of $\widetilde{V}$ into  itself such that 
$\widetilde{\gamma }(\widetilde{A}\widetilde{u}, 
\widetilde{A}\widetilde{v}) = 
\widetilde{\gamma }(\widetilde{u},\widetilde{v})$ for every 
$\widetilde{u},\widetilde{v} \in \widetilde{V}$. 
$\OO (\widetilde{V}, \widetilde{\gamma })$ is called 
the \emph{real orthogonal group}. 
It is a Lie group whose Lie algebra $\oo (\widetilde{V}, \widetilde{\gamma})$ is the set of real linear maps $\widetilde{\xi}:\widetilde{V} \rightarrow 
\widetilde{V}$ such 
that $\widetilde{\gamma }(\widetilde{\xi}\widetilde{u},
\widetilde{v}) + \widetilde{\gamma }(\widetilde{u}, 
\widetilde{\xi}\widetilde{v}) = 0$ for every $\widetilde{u},
\widetilde{v} \in \widetilde{V}$. Suppose that $\gamma $ is a symmetric bilinear form on a real finite dimensional vector space $V$, 
whose kernel is the span of a nonzero vector $v^0$. 
Let ${\OO (V, \gamma )}_{v^0}$ be the set of 
bijective real linear mapings $A:V \rightarrow V$ such that 
$Av^0 = v^0$ and $\gamma (Au, Av) = \gamma (u,v)$ for every $u,v \in V$. 
${\OO (V, \gamma )}_{v^0}$ is called the \emph{generalized real orthogonal 
group}. It is a Lie group whose Lie algebra 
${\oo (V, \gamma )}_{v^0}$ consists of real linear maps $\xi :V \rightarrow V$ 
such that $\xi v^0 =0$ and $\gamma (\xi u, v) + \gamma (u, \xi v) =0$ 
for every $u,v \in V$. \medskip 

\section{Basic concepts}

The goal of this paper is to find a unqiue 
representative (= normal form) for each orbit of the adjoint action of 
${\OO (V, \gamma )}_{v^0}$ on ${\oo (V, \gamma )}_{v^0}$. \medskip 

We begin by defining some basic concepts, which follow 
\cite{burgoyne-cushman} and \cite{cushman-vanderkallen}. A 
\emph{pair} $(\widetilde{\xi }|\widetilde{W} , 
\widetilde{W}; \widetilde{\gamma }|\widetilde{W})$ 
is a proper $\widetilde{\xi }$-invariant 
subspace $\widetilde{W}$ of $\widetilde{V}$ 
such that $\widetilde{\gamma }$ is nondegenerate on $\widetilde{W}$ 
and $\widetilde{\xi }\in \oo (\widetilde{W}, 
\widetilde{\gamma }|\widetilde{W})$. We 
say that two pairs $(\widetilde{\xi }|\widetilde{W}, 
\widetilde{W}; \widetilde{\gamma }|\widetilde{W})$ and 
$(\widetilde{\xi }'|{\widetilde{W}}', {\widetilde{W}}'; 
{\widetilde{\gamma }}'|{\widetilde{W}}' )$ are \emph{equivalent} if there is a bijective real linear mapping $\widetilde{P}:\widetilde{W} \rightarrow 
{\widetilde{W}}'$ such that 
${\widetilde{\gamma }}'(\widetilde{P}\widetilde{u}, 
\widetilde{P}\widetilde{v}) = 
\widetilde{\gamma }(\widetilde{u},\widetilde{v})$  and 
$(\widetilde{P}\comp \widetilde{\xi } )\widetilde{w} = 
({\widetilde{\xi }}' \comp \widetilde{P})\widetilde{w}$ for every 
$\widetilde{u},\widetilde{v},\widetilde{w} \in \widetilde{W}$. Being 
equivalent is an equivalence relation on the set of pairs. An 
equivalence class $\Delta $ of pairs is called a \emph{type}. If 
${\widetilde{W}}_i$, $i=1,2$, are proper $\widetilde{\xi }$-invariant, and 
$\widetilde{\gamma }$-orthogonal 
subspaces of $\widetilde{W}$, whose direct sum is $\widetilde{W}$ 
and on which 
$\widetilde{\gamma }|{\widetilde{W}}_i $ is nondegenerate, then the type $\Delta $, represented by $(\widetilde{\xi }|\widetilde{W},
\widetilde{W}; \widetilde{\gamma }|\widetilde{W})$, 
is the \emph{sum} of two types ${\Delta }_1$ and ${\Delta }_2$ represented by 
$(\widetilde{\xi }|{\widetilde{W}}_1, {\widetilde{W}}_1; 
\widetilde{\gamma }|{\widetilde{W}}_1)$ and 
$(\widetilde{\xi }|{\widetilde{W}}_2, {\widetilde{W}}_2; 
\widetilde{\gamma }|{\widetilde{W}}_2)$, respectively. 
We write $\Delta = {\Delta }_1 + {\Delta }_2$. A type, which 
cannot be written as the sum of two types, is \emph{indecomposable}. Indecomposable types have been classified in \cite{burgoyne-cushman}.  \medskip 

A \emph{triple} $(W, \xi |W , v^0; \gamma |W)$ is a proper, $\xi $-invariant 
subspace $W$ of $V$ containing $v^0$ such that $\ker \gamma |W = 
\spann \{ v^0 \}$ and $\xi |W \in {\oo (W, \gamma |W )}_{v^0}$. Two triples 
$(W, \xi |W, v^0; \gamma |W)$ and 
$(W', {\xi }'|W' ,  (v')^0; {\gamma }'|W')$ are \emph{equivalent} 
if there is a bijective real linear mapping $P:W \rightarrow W'$ 
such that $Pv^0 = (v')^0$, ${\gamma }'(Pu, Pv) = 
\gamma (u,v)$ and $({\xi }' \comp P)w = 
(P \comp \xi )w$ for every $u,v,w\in W$. Being equivalent is an 
equivalence relation on the set of triples. An equivalence class 
$\underline{\underline{\Delta }}$ of triples is called a 
\emph{special type}. If $(W, \xi |W , v^0; \gamma |W)$ is a triple, 
representing $\underline{\underline{\Delta}}$, with 
$\xi $ nilpotent, then $\underline{\underline{\Delta }}$ 
is a \emph{nilpotent special type}. 
If $W_i$, $i=1,2$, are proper $\xi $-invariant, 
$\gamma $-orthogonal subspaces of $W$ with $v^0 \in W_1$, $W = 
W_1\oplus W_2$ and $\ker \gamma |W_1= \spann \{ v^0 \} $, then the special type 
$\underline{\underline{\Delta }}$, represented by $(W, \xi |W, v^0; \gamma |W)$, 
is the \emph{sum} of the special 
type ${\underline{\underline{\Delta }}}'$, represented by 
$(W_1, \xi |W_1, v^0; \gamma |W_1)$, and the type $\Delta $ represented by 
$(\xi |W_2, W_2; \gamma |W_2)$. We write $\underline{\underline{\Delta }} = 
{\underline{\underline{\Delta}}}' + \Delta$. A special type is 
\emph{indecomposable} if it cannot be written as the sum of 
a special type and a type. \medskip

For $\xi \in {\oo (V, \gamma )}_{v^0}$ let $V_0$ be the direct sum of 
generalized eigenspaces
of $\xi $ each corresponding to an eigenvalue $0$ and let $U$ be the 
direct sum of generalized eigenspaces of $\xi $ each corresponding to a 
nonzero eigenvalue of $\xi $. Then $V_0$ and $U$ are $\xi $-invariant, 
$\gamma $-orthogonal subspaces of $V$ such that $V = V_0 \oplus U$. 
Here $v^0 \in V_0$, since $\xi v^0 =0$, and $\ker \gamma |V_0 = 
\spann \{ v^0 \}$.  Thus $\gamma |U$ is nondegenerate. Therefore the special type 
$\underline{\underline{\Delta }}$, 
represented by $(V, \xi , v^0; \gamma )$, is the sum of the special 
type ${\underline{\underline{{\Delta }}}}_0$, represented by 
$(V_0, \xi |V_0 , v^0; \gamma |V_0)$, and the type $\Delta $, 
represented by the pair $(\xi |U, U; \gamma |U)$, that is, 
$\underline{\underline{\Delta}} = {\underline{\underline{\Delta }}}_0 + 
\Delta $. Since $\xi |V_0$ is nilpotent, ${\underline{\underline{\Delta }}}_0$ is a nilpotent special type. \medskip

So we have proved \medskip 

\noindent \textbf{Proposition 1.} Every special type $\underline{\underline{\Delta}}$ 
may be written uniquely as the sum of a nilpotent special type 
${\underline{\underline{\Delta}}}_0$ and a type $\Delta $. \medskip  

Let $(V, \xi , v^0; {\gamma })$ be a triple representing the special type 
$\underline{\underline{\Delta }}$. Since $v^0 \in \ker \xi $, 
the vector $v^0$ lies at the end of a Jordan 
chain in $V$. Let $h$ be the largest nonegative integer such that there 
is a vector $v \in V$ with ${\xi}^hv = v^0$. Then a longest Jordan 
chain in $V$ ending at $v^0$ has length $h+1$, since ${\xi }^{h+1}v= 
\xi v^0 = 0$. We call $h$ the \emph{special height} of 
$\underline{\underline{\Delta }}$ 
and denote it by $\mathrm{sht}\, \underline{\underline{\Delta }}$. 
The special height of $\underline{\underline{\Delta }}$ does 
not depend on the choice of triple representing the special type 
$\underline{\underline{\Delta }}$. \medskip 

\noindent \textbf{Lemma 2.} Let $\underline{\underline{\Delta }}$ be a special 
type with special height $h$. Suppose that 
\begin{equation}
\underline{\underline{\Delta }} = 
\underline{\underline{\Delta }}' + \Delta , 
\label{eq-sec1one}
\end{equation}
where $\underline{\underline{\Delta}}'$ is a special nilpotent type 
and $\Delta $ is a type. Then 
$\mathrm{sht}\, \underline{\underline{\Delta }}'$ $= 
\mathrm{sht}\, \underline{\underline{\Delta }}$ and $\mathrm{ht}\, \Delta < 
\mathrm{sht}\, {\underline{\underline{\Delta }}}'$. 
Here $\mathrm{ht}\, \Delta $ is the height of the type $\Delta $, see 
\cite{burgoyne-cushman}. \medskip  

\noindent \textbf{Proof.} Suppose that the triple $(V,\xi , v^0; {\gamma })$ 
represents the special type $\underline{\underline{\Delta }}$. 
Since $\mathrm{sht}\, \underline{\underline{\Delta }} =h$, there is a vector 
$v \in V$ such that ${\xi }^{h}v = v^0$. 
Let $(V_1, \xi |V_1, v^0;$ $\gamma |V_1)$ 
be a triple representing the special type $\underline{\underline{\Delta }}'$ and $(\xi |V_2, V_2;\gamma |V_2)$ be a pair representing the type $\Delta $. Then $V = V_1 \oplus V_2$, $V_i$ are $\xi $-invariant, and $v^0 \in V_1$. 
Write $v = v_1 +v_2 \in V_1 \oplus V_2$. Then $v^0 = {\xi }^{h}v 
= {\xi }^{h}v_1 + {\xi }^{h}v_2 \in V_1 \oplus V_2$. 
But $v^0 \in V_1$. So ${\xi }^{h}v_2=0$, 
which implies that the height $\mathrm{ht}\, \Delta $ of $\Delta $ is 
less than $h$. Clearly $\mathrm{sht}\, {\Delta }' 
\le \mathrm{sht}\, \Delta = h$. But $v^0 = {\xi }^{h}v_1$. This 
implies $\mathrm{sht}\, \underline{\underline{\Delta '}} 
\ge h$. \hfill $\square $ \medskip 

Using lemma 2 repeatedly we may assume that 
$\underline{\underline{\Delta }}'$ in (\ref{eq-sec1one}) is an 
\emph{indecomposable} nilpotent special type. \medskip 

Suppose that $\underline{\underline{\Delta }}$ is a special type, 
represented by $(V, \xi , v^0; \gamma )$, which is 
\emph{not} equal to the nilpotent special type $\underline{\underline{\tau}}$, represented by $(\spann \{ v^0 \} , 0 , v^0; 0)$.  Let 
$\overline{V} = V/\spann \{ v^0 \} $ with projection map 
$\pi : V \rightarrow \overline{V}:v \mapsto \overline{v}$. Since 
$\underline{\underline{\Delta}} \ne \underline{\underline{\tau}}$, 
it follows that $\overline{V} \ne \{ \overline{0} \} $. 
Let $\overline{\gamma }$ be the bilinear form on $\overline{V}$ defined by 
$\overline{\gamma }(\overline{v}, \overline{w}) = \gamma (v,w)$, 
where $v\in \overline{v}\in \overline{V}$ and $w \in \overline{w} \in \overline{V}$. The form 
$\overline{\gamma }$ is well defined and is symmetric. In fact it is nondegenerate. For if $\overline{\gamma }(\overline{v}, \overline{w}) =0$ for 
every $\overline{w} \in \overline{V}$, then $0 = \gamma (v,w)$ for 
every $w \in V$. Therefore, $v \in \ker \gamma  = \spann \{ v^0 \} $, 
that is, $v = \mu \, v^0$ for some $\mu \in \R $. So $\overline{v} 
= \overline{0}$. Since $\xi v^0 =0$, the map $\xi \in {\oo (V, \gamma )}_{v^0}$ 
induces a linear map $\overline{\xi }:\overline{V} \rightarrow 
\overline{V}: \overline{v} \mapsto \overline{\xi v}$. Because 
\begin{displaymath}
\overline{\gamma }(\overline{\xi }\overline{v}, \overline{w}) 
= \gamma (\xi v, w) = 
-\gamma (v, \xi w) 
= - \overline{\gamma }(\overline{v}, \overline{\xi }\overline{w}),  
\end{displaymath}
for every $\overline{v}, \overline{w} \in \overline{V}$, 
it follows that $\overline{\xi } \in 
\oo (\overline{V}, \overline{\gamma })$. Therefore $(\overline{\xi }, 
\overline{V}; \overline{\gamma })$ represents a type 
$\overline{\Delta }$, called the \emph{type induced} by the special type 
$\underline{\underline{\Delta }}$. \medskip 

\noindent \textbf{Lemma 3.} Suppose that $\underline{\underline{\Delta }}$ is a 
special type, represented by $(V,\xi ,$ $v^0;\gamma )$, 
with special height $h+1$. Then the induced type 
$\overline{\Delta }$ has height $h$. \medskip 

\noindent \textbf{Proof.} Because 
$\mathrm{sht}\, \underline{\underline{\Delta }}$ is 
$h+1$, there is a vector $v \in V$ such that ${\xi }^{h+1}v =v^0$. Since 
${\xi }^h v \not \in \spann \{ v^0 \}$, it follows that 
${\overline{\xi }}^h \overline{v} \ne 
\overline{0}$. But ${\overline{\xi}}^{h+1}\overline{v} = 
\overline{{\xi}^{h+1}v} = \overline{v^0} =\overline{0}$. Thus the induced type 
$\overline{\Delta }$, represented by $(\overline{\xi }, \overline{V}; \overline{\gamma })$, 
has height $h$. \hfill $\square $    

\section{Classification of indecomposable nilpotent \\ special types}

In this section we classify indecomposable nilpotent special types. Our 
argument shows that an indecomposable nilpotent special type is 
uniquely determined by its induced type.  
 
\subsection{A rough classfication} 

In this subsection we give a rough classification of 
indecomposable nilpotent special types. \medskip 

\noindent \textbf{1}. We begin by treating the case when the indecomposable 
nilpotent special type $\underline{\underline{\Delta }}$, represented by $(V,N,v^0; \gamma )$, has special height $0$. By results of \cite{burgoyne-cushman}, the underlying $\gl $-type 
$\Delta $, represented by the pair $(V,N)$, is the unique sum of indecomposable $\gl $-types ${\Delta }_0 +{\Delta }_1 + \cdots + {\Delta }_r$ with heights $h_0 \le \cdots \le h_r$, respectively. For $0 \le i \le r$ let $(V_i, N|V_i)$ be a pair representing the 
$\gl$-type ${\Delta }_i$. Since the special height of $\underline{\underline{\Delta }}$ is 
$0$, the longest Jordan chain in $V$ which ends at $v^0$ is 
$\{ v^0 \}$. Clearly the pair $(V_0=\spann \{ v^0 \},0)$ represents  
an indecomposable $\gl $-type of height $0$. Therefore $h_0 =0$ 
and $(\spann \{ v^0 \},0) \in {\Delta }_0$. Write $V = V_0 \oplus U$, where 
$U = V_1 \oplus \cdots \oplus V_r$. Then $V_0$ and $U$ are $N$-invariant 
with $v^0 \in V_0$. Since $\ker \gamma = \spann \{ v^0 \}$, it follows 
that $V_0$ and $U$ are $\gamma $-orthogonal and $\gamma |U$ is 
nondegenerate. Therefore $\underline{\underline{\Delta}}$ is the 
sum of the nilpotent special type 
$\underline{\underline{\tau }}$, represented by 
$(\spann \{ v^0 \}, 0, v^0; 0)$, and a type represented by 
$(N|U, U; \gamma |U)$. But by hypothesis 
$\underline{\underline{\Delta}}$ is indecomposable. Therefore 
$\underline{\underline{\Delta }} = \underline{\underline{\tau }}$. This 
completes alternative 1. \medskip 

\noindent \textbf{2}. Suppose that $(V, N, v^0; \gamma )$ 
represents an indecomposable special nilpotent type 
$\underline{\underline{\Delta }}$ of special height $h+1$. Let $\{ w,\, Nw, \, \ldots , N^{h+1}w =v^0 \}$ be a longest Jordan chain ending at $v^0$ which lies in $V$. There are two possibilities. \medskip 

\noindent \textbf{a}. $\gamma (w, N^hw ) \ne 0$. This hypothesis implies that $h$ is even. If not, then $h$ is odd and 
\begin{displaymath}
\gamma (w, N^hw) = (-1)^h \, \gamma (N^h w, w) = -\gamma (w, N^hw).
\end{displaymath}
Therefore $\gamma (w, N^hw) =0$, which contradicts our hypothesis. Since 
$N^h w \not \in \spann \{v^0 \}$, it follows that 
${\overline{N}}^h\overline{w} \ne 0$. Also  
${\overline{N}}^{h+1}\overline{w} =0$, because $N^{h+1}w =v^0$.  So 
$\{ \overline{w}, \, \overline{N}\overline{w}, \, \ldots , 
\, {\overline{N}}^h\overline{w} \}$ is a 
Jordan chain in $\overline{V}$ of length $h+1$. Let $\overline{U}$ be the space spanned by the elements of this chain. Let $\overline{G}$ be the Gram matrix 
of $\overline{\gamma }$ on $\overline{U}$. If $i+j \ge h+1$, then 
\begin{displaymath}
{\overline{\gamma }}({\overline{N}}^i\overline{w}, {\overline{N}}^j
\overline{w}) = (-1)^i\, \overline{\gamma }(\overline{w}, 
{\overline{N}}^{i+j}\overline{w}) =0,
\end{displaymath}
since ${\overline{N}}^{h+1}\overline{w} =0$. Therefore 
all the entries of $\overline{G}$ below 
the antidiagonal are equal to $0$. All the entries of $\overline{G}$ on 
the antidiagonal are nonzero, because   
\begin{displaymath}
\overline{\gamma }({\overline{N}}^i\overline{w}, 
{\overline{N}}^{h-i}\overline{w}) = (-1)^i \, \overline{\gamma } ( 
\overline{w}, {\overline{N}}^h\overline{w}) \ne 0, \quad 
\mbox{for $0 \le i \le h$,}
\end{displaymath} 
since $\overline{\gamma }(\overline{w}, {\overline{N}}^h\overline{w}) = \gamma (w, N^hw) \ne 0$ by hypothesis. Therefore $\det \overline{G} \ne 0$, that is, $\overline{\gamma }|\overline{U}$ is 
nondegenerate. Since $\overline{U}$ is spanned by a Jordan chain 
of length $h+1$, it follows that $\overline{\Delta}$ is a uniform type of height 
$h$, that is, $\ker {\overline{N}}^h = \overline{N} \, \overline{U}= 
\overline{NU}$. Therefore $\overline{\gamma }$ induces a nondegenerate 
bilinear form $\widehat{\gamma }$ on $\widehat{U} = \overline{U}/ 
\overline{NU}$ defined by $\widehat{\gamma }(\widehat{z}, 
\widehat{w}) = \overline{\gamma }(\overline{z}, {\overline{N}}^h 
\overline{w})$, where $\widehat{z} \in \overline{z}$, $\widehat{w} 
\in \overline{w}$. The bilinear form $\widehat{\gamma}$ is symmetric, since $h$ is even.  
Because $\dim \widehat{U} =1$, we can choose a vector $\widehat{w} 
\in \widehat{U}$ so that $\widehat{\gamma }(\widehat{w}, 
\widehat{w}) = \varepsilon $ with ${\varepsilon }^2 =1$. In other words, 
with respect to the basis $\{ \widehat{w} \} $ of $\widehat{U}$ the Gram matrix of 
$\widehat{\gamma }$ is $(\varepsilon )$. From the classification 
of indecomposable types in \cite{burgoyne-cushman} we see that the induced type 
$\overline{\Delta }$, represented by $(\overline{N}|\overline{U}, \overline{U}; \overline{\gamma }|\overline{U})$, is ${\Delta }^{\varepsilon }_h(0)$. This completes case \textbf{a}. \medskip 

\noindent \textbf{b}. $\gamma (w, N^hw) =0$. Since 
$N^hw \not \in \spann \{ v^0 \}$, it follows that 
${\overline{N}}^h\overline{w} \ne \overline{0}$. Moreover, 
${\overline{N}}^{h+1}\overline{w} = \overline{N^{h+1}w} = \overline{v^0} 
= \overline{0}$. Therefore  $\{ \overline{w}, \, \overline{N}\overline{w}, 
\, \ldots ,$ ${\overline{N}}^h\overline{w} \} $ is a Jordan chain in 
$\overline{V}$. Because $\overline{\gamma }$ is nondegenerate on 
$\overline{V}$, there is a nonzero vector $\overline{z} \in \overline{V}$ 
such that $\overline{\gamma }(\overline{z}, {\overline{N}}^h\overline{w}) 
\ne 0$. Therefore ${\overline{N}}^h \overline{z} \ne \overline{0}$, for  
otherwise $\overline{0} = \overline{\gamma }({\overline{N}}^h \overline{z}, 
\overline{w}) = (-1)^h \overline{\gamma }(\overline{z}, 
{\overline{N}}^h\overline{w})$, 
which is a contradiction. Let $\overline{U} = \spann \{ \overline{z}, \, 
\ldots , {\overline{N}}^h\overline{z}, \, \overline{w}, \ldots , 
{\overline{N}}^h\overline{w} \} \subseteq \overline{V}$. Since 
$\mathrm{ht}\, \overline{\Delta } =h$, 
it follows that ${\overline{N}}^{h+1} 
\overline{z} =0$. Therefore $\overline{U}$ is $\overline{N}$-invariant. 
Consider the basis $ \{ \overline{z}, 
\overline{w}, \overline{N}\overline{z}, \overline{N}\overline{w}, 
\ldots , {\overline{N}}^h \overline{z} , {\overline{N}}^h\overline{w} \}$ 
of $\overline{U}$. Since 
\begin{eqnarray*}
\gamma ({\overline{N}}^i\overline{w}, {\overline{N}}^j \overline{w}) 
& = & 0, \quad \mbox{if $i+j \ge h+1$} \\ 
\gamma ({\overline{N}}^i\overline{z}, {\overline{N}}^j \overline{z}) 
& = & 0, \quad \mbox{if $i+j \ge h+1$} \\ 
\gamma ({\overline{N}}^i\overline{z}, {\overline{N}}^j \overline{w}) 
& = & 0, \quad \mbox{if $i+j \ge h+1$} \\
\gamma ({\overline{N}}^i\overline{w}, {\overline{N}}^j \overline{z}) 
& \ne  & 0, \quad \mbox{if $i+j = h$},
\end{eqnarray*}
the Gram matrix $\overline{G}$ of $\overline{\gamma }|\overline{U}$ has zero
entries below the antidiagonal and nonzero entries on the antidiagonal. Therefore $\overline{\gamma }$ is nondegenerate on $\overline{U}$. 
Consequently the pair $(\overline{N}|\overline{U}, \overline{U}; 
\overline{\gamma }|\overline{U})$ represents a type 
$\overline{\Delta }$, which is nilpotent of height $h$. 
Since $\overline{U}$ is spanned by two Jordan chains of length $h+1$, 
it follows that $\overline{\Delta }$ is a uniform type of height $h$, that is, $\ker {\overline{N}}^h = 
\overline{N}\, \overline{U}= \overline{NU}$. Therefore $\overline{\gamma }$ induces a nondegenerate bilinear form $\widehat{\gamma }$ on 
$\widehat{U} = \overline{U}/\overline{N}\overline{U}$ 
defined by $\widehat{\gamma }(\widehat{z}, \widehat{w}) = 
\overline{\gamma }(\overline{z}, {\overline{N}}^h\overline{w})$, 
where $\widehat{z} \in \overline{z} \in \overline{V}$, $\widehat{w} \in \overline{w} \in 
\overline{V}$. The bilinear form $\widehat{\gamma }$ is symmetric if $h$ is even and skew symmetric if $h$ is odd. Since $\dim \widehat{U} =2$, we can choose vectors $\widehat{z}, \widehat{w} \in \widehat{U}$ 
so that the Gram matrix of $\widehat{\gamma }$ is 
{\tiny $\begin{pmatrix} 0 & 1 \\ 1 & 0 \end{pmatrix}$} if $h$ is even and 
{\tiny $\begin{pmatrix} 0 & -1 \\ 1 & 0 \end{pmatrix}$} if $h$ is odd. Thus $\widehat{\gamma }(\widehat{z}, 
\widehat{w}) =1$ in both cases. From the classification of indecomposable types in \cite{burgoyne-cushman} $\overline{\Delta }$ is ${\Delta }^{+}_h(0) + {\Delta }^{-}_h(0)$ if $h$ is even and 
${\Delta }_h(0,0)$ if $h$ is odd. This completes case \textbf{b}. \medskip 

The rough classification of indecomposable nilpotent special 
types is now established.   \hfill $\square $ \medskip

\subsection{Fine classfication}

In this subsection we give a fine classification of indecomposable 
nilpotent special types.  \medskip 

\noindent \textbf{Theorem 5.} An indecomposable nilpotent special  
type $\underline{\underline{\Delta}}$, represented by the triple 
$(V,N, v^0; \gamma )$, is exactly one of the following. 
\begin{itemize}
\item[1.] $\underline{\underline{\tau}}$. Here 
$\mathrm{sht}\, \underline{\underline{\tau}} = 0$.   
Then $V =\spann \{ v^0 \} $. With respect to this basis the 
matrix of $\gamma $ is $0$ and the matrix of $N$ is zero.  
\item[2.] ${\underline{\underline{{\Delta }^{\varepsilon}_{h+1}}}}(0)$, 
$\alpha >0$. Here $h$ is even, $\varepsilon = \pm$, and $\mathrm{sht}\, {\underline{\underline{{\Delta }^{\varepsilon}_{h+1}}}}(0) = h+1$. We have 
$v^0 = \alpha N^{h+1}w$ with $\alpha >0$. We call $\alpha $ a 
\emph{modulus}. If $h=0$, then 
$\mathfrak{f} = \{ \alpha Nw; \, w \} $ is a basis of $V$. With respect 
to $\mathfrak{f}$ the Gram matrix of $\gamma $ is {\tiny $ 
\left( \begin{array}{c|c} 0 & 0 \\ \hline 
0 & \varepsilon \, 1 \end{array} \right) $} and the matrix of $N$ is 
{\tiny $ \left( \begin{array}{c|c} 0 & {\alpha }^{-1} \\ \hline 
0 & 0 \end{array} \right) $}. When $h$ is nonzero    
\begin{align}
\mathfrak{f} & = \mbox{$\{ \alpha N^{h+1}w;N^{h/2-1}w, \, N^{h/2-2}w, \, \ldots , w; \, 
N^{h/2}w; $} \notag \\
&\hspace{.5in}\mbox{$ \delta \, N^{h/2+1}w, 
-\delta \, N^{h/2+2}w, \ldots , \, \delta N^hw \} $} \notag   
\end{align}
is a basis of $V$ with $\delta = \varepsilon (-1)^{h/2}$, where $\gamma (w, N^hw) = 
\varepsilon = \pm 1$. With respect to the basis $\mathfrak{f}$ the Gram matrix of $\gamma $ is  
\begin{equation}
\left( \begin{array}{c|c}
0 & \\ \hline 
  & {\mathcal{I}}_{h, \delta } 
\end{array} \right), \, \, \, \mathrm{where} \, \, \, 
{\mathcal{I}}_{h, \delta } = 
\mbox{{\tiny $\left( \begin{array}{c|c|c}
   0     &    0     & I_{h/2} \\ \hline 
   0     & \delta   &  0      \\  \hline 
    I_{h/2} &    0     &  0
\end{array} \right) $.}} 
\label{eq-thm5pt2gram}
\end{equation}
The index of $\gamma $ is $h/2$ if $\delta =1$ and $h/2+1$ if $\delta =-1$. 
The matrix of $N$ with respect to the basis $\mathfrak{f}$ is 
\begin{equation}
\mbox{{\footnotesize $\left( \begin{array}{c|c|c|r}
0 &       &                 & \delta {\alpha }^{-1} e^T_{h/2} \\ \hline 
  &J_{h/2}&                 &                              \\ \hline
  &\rule[-2pt]{0pt}{12pt}e^T_1    &                 &            \\ \hline
  &                         &-\delta e_{1} 
  &\rule[-2pt]{0pt}{12pt} -J^T_{h/2}                   
\end{array} \right) .$}}
\label{eq-thm5pt2nil}
\end{equation}
Here $J_{h/2 }$ is the $h/2 \times h/2$ upper Jordan block.  
\item[3.] $\underline{\underline{{\Delta }_{h+1}(0,0)}}$ if $h$ is odd, or 
$\underline{\underline{{\Delta }^{+}_{h+1}(0) + {\Delta }^{-}_{h+1}(0)}}$ 
if $h$ is even. Here \linebreak 
$\mathrm{sht}\, \underline{\underline{{\Delta }_{h+1}(0,0)}} = 
\mathrm{sht}(\underline{\underline{{\Delta }^{+}_{h+1}(0) + {\Delta }^{-}_{h+1}(0)}}) 
= h+1$. We have $v^0 = N^{h+1}v$.   
\begin{displaymath}
\mathfrak{f} = \mbox{$\{ N^{h+1}v;\,  v, \, Nv, \, 
\ldots N^h v;\, N^hw, -N^{h-1}w, \, \ldots , (-1)^hw \} $,}
\end{displaymath}
is a basis of $V$. The Gram matrix of $\gamma $ 
with respect to $\mathfrak{f}$ is 
\begin{equation}
\mbox{{\footnotesize $\left( \begin{array}{c|c|c} 
0 & 0 & 0 \\ \hline 
0 & 0 & I_{h+1} \\ \hline
0 & I_{h+1} & 0 \end{array} \right) $}} 
\label{eq-thm5pt3gram}
\end{equation}
and the matrix of $N$ is 
\begin{equation}
\mbox{{\footnotesize $\left( \begin{array}{c|l|c}
0 & e^T_{h+1} &  \\ \hline 
\rule{0pt}{10pt}  &  J^T_{h+1} & \\ \hline
\rule{0pt}{10pt}  &          & -J_{h+1} \end{array} \right) $}}, 
\label{eq-thm5pt3nil}
\end{equation}    
where $J_{h+1}$ is an $(h+1) \times (h+1)$ upper Jordan block.  
\end{itemize} 

\noindent \textbf{Proof.} Suppose that $(V,N,v^{0}; \gamma )$ represents 
an indecomposable nilpotent special type $\underline{\underline{\Delta }}$, whose induced type $\overline{\Delta }$ is given. \medskip 

\noindent 1. The case when $\underline{\underline{\Delta }}$ is   
$\underline{\underline{\tau }}$ has been dealt with in the first 
case of the rough classification in \S 3.1. This establishes alternative 1. \medskip 

\noindent 2. By case \textbf{a} of the rough classification, 
the induced type $\overline{\Delta }$ of height $h$,  where $h$ is even, 
and is represented by the pair 
$(\overline{N}, \overline{V};\overline{\gamma })$. We can choose 
$\overline{w} \in \overline{V}$ so that $\overline{V}$ has a basis 
\begin{equation}
\mbox{{\footnotesize $ \{ {\overline{N}}^{h/2-1}\overline{w}, 
{\overline{N}}^{h/2-2}\overline{w},\ldots , \overline{w}; 
{\overline{N}}^{h/2}\overline{w}; -\delta {\overline{N}}^{h/2+1}\overline{w}, 
\, \delta {\overline{N}}^{h/2+2}\overline{w}; \, \ldots ,\, 
\delta {\overline{N}}^h\overline{w} \} $.}}  
\label{eq-thm5pt2one}
\end{equation} 
Here $\overline{\gamma }({\overline{N}}^i\overline{w} ,{\overline{N}}^{h-i}\overline{w}) = \varepsilon (-1)^i$ for $0 \le i \le h$ and the value of 
$\overline{\gamma }$ is $0$ on any other pair of vectors in 
the basis (\ref{eq-thm5pt2nil}). Also 
$\delta = \varepsilon (-1)^{h/2}$.  \medskip 

We now reconstruct the special type $\underline{\underline{\Delta}}$ 
from its induced type $\overline{\Delta}$. 
Since $\mathrm{sht}\, \underline{\underline{\Delta }} = h+1$, 
the vector space $V$ has a basis $\mathfrak{f}=${\footnotesize 
$ \{w, \, Nw, \, \ldots , N^hw, \, N^{h+1}w  \} $}, where 
$w \in \overline{w} \in \overline{V}$. Therefore $V$ is spanned by a Jordan chain 
of length $h+2$. Because ${\overline{N}}^{h+1} \overline{w} =0$, 
there is a real number $\mu $ such that $N^{h+1}w = \mu  \, v^0$. 
Moreover, we can choose $w \in \overline{w}$ so that 
$\gamma (N^iw, N^{h-i}w) = (-1)^i \varepsilon$ for $0 \le i \le h$ and the value of $\gamma $ is $0$ on any other pair of vectors in the basis $\mathfrak{f}$. If $\mu =0$, 
then the triple $(V, N,v^0;\gamma )$ does not represent a nilpotent special type of special height $h+1$ because there is no Jordan chain 
in $V$ of length $h+2$. Therefore $\mu  \ne 0$. Let $\alpha = {\mu }^{-1}$. With respect to the basis 
\begin{equation}
\begin{array}{l}
\{ \alpha N^{h+1}w; \, N^{(h-1)/2}w, \, 
N^{(h-2)/2} w,\ldots , w; N^{h/2} w;  \\
\hspace{1.5in}\delta \, N^{h/2+1} w, \, -\delta N^{h/2+2}w, \ldots , \delta N^hw \} ,
\end{array} 
\label{eq-thm5pt2two}
\end{equation}
the Gram matrix of $\gamma $ is (\ref{eq-thm5pt2gram}). Changing the sign 
of the vector $w$ does not change the Gram matrix of $\gamma $. Therefore 
we may select the modulus $\alpha $ so that it is positive. It is 
straightforward to see that the matrix of $N$ with respect to the basis 
(\ref{eq-thm5pt2two}) is (\ref{eq-thm5pt2nil}). \medskip 

\noindent 3. By case \textbf{b} of the rough classification, the induced type 
$\overline{\Delta }$ has height $h$ and is represented by the pair 
$(\overline{N}, \overline{V};\overline{\gamma })$. 
We now reconstruct the special type $\underline{\underline{\Delta}}$ from its induced type 
$\overline{\Delta}$. From case 3 of the rough classification in \S 3.1 we can choose 
$\overline{z} \in \widehat{z} \in \overline{V}$ and $\overline{w} \in \widehat{w} \in 
\overline{V}$ so that  
\begin{equation}
\{ \overline{z},\, \overline{N}\overline{z}, \, \ldots , 
{\overline{N}}^h\overline{z};\, 
\overline{w},\, \overline{N}\overline{w}, \, \ldots , 
{\overline{N}}^h\overline{w} \}
\label{eq-thm5pt3onestar}
\end{equation}
is a basis of $\overline{V}$ such that 
$\overline{\gamma }({\overline{N}}^i\overline{z}, 
{\overline{N}}^{h-i}\overline{w}) = (-1)^i$ for $0 \le i \le h$ 
and the value of $\overline{\gamma }$ on 
any other pair of vectors in the basis (\ref{eq-thm5pt3onestar}) 
is zero. \medskip

Let $z \in \overline{z} \in \overline{V}$ and $w \in \overline{w}\in \overline{V}$. Since 
${\overline{N}}^{h+1}\overline{z} = 
{\overline{N}}^{h+1}\overline{w} =\overline{0}$, 
there are real numbers $\lambda $ and 
$\mu $ such that $N^{h+1}z = \lambda \, v^0$ and $N^{h+1}w = \mu \, v^0$. 
The vectors 
\begin{equation}
\mbox{{\footnotesize $ \{ z,\, Nz, \, \ldots , N^{h+1}z = \lambda \, v^0; 
w, \, Nw, \, \ldots , N^hw, \, N^{h+1}w = \mu \, v^0 \} $ }}
\label{eq-thm5pt3two}
\end{equation} 
span $V$. Suppose that $\lambda = \mu =0$ in (\ref{eq-thm5pt3two}). Then there is no 
Jordan chain in $V$ of length $h+2$. This contradicts the fact 
that the special height of $\underline{\underline{\Delta }}$ is $h+1$. Therefore not both $\lambda $ and $\mu $ are zero. Interchanging 
$z$ and $w$ in (\ref{eq-thm5pt3two}), if 
necessary, we may suppose that $\mu \ne 0$. Let $\eta = \mu z - \lambda w$ and 
$\zeta = {\mu }^{-1}w$. Then $N^{h+1}\eta =0$ and $N^{h+1}\zeta = v^0$. 
Note that 
\begin{equation}
\{ \eta , \, N\eta , \, \ldots , N^h\eta ; \, \zeta, \, N\zeta , \, \ldots , 
N^{h+1}\zeta = v^0 \} 
\label{eq-thm5pt3three}
\end{equation} 
is a basis of $V$. Therefore $V$ is spanned by two Jordan chains one 
of length $h+2$, which ends at $v^0$, and the other of length $h+1$. \medskip

We now calculate the Gram matrix of $\gamma $ with 
respect to the basis (\ref{eq-thm5pt3three}). When $i+j = h$ we have 
\begin{displaymath}
\gamma (N^i\eta , N^j\zeta ) = \gamma (\mu N^iz - \lambda N^i w, 
{\mu }^{-1}N^jw) 
\, =  (-1)^i \gamma (z, N^hw ) \, = \, (-1)^i;
\end{displaymath}
while when $i+j \ne h$ we obtain 
\begin{eqnarray*}
\gamma (N^i\eta , N^j \zeta ) & = & (-1)^i \gamma (z, N^{i+j}w) \, 
= \, 0 ; \\
\gamma (N^i\eta , N^j \eta ) & = & (-1)^i \gamma (\eta , N^{i+j}\eta ) \\ 
& = & (-1)^i[ {\mu }^2\, \gamma (z, N^{i+j}z) - 2\lambda \mu \, 
\gamma (z, N^{i+j}w) 
+{\lambda }^2\, \gamma (w,N^{i+j}w)] \\
& = & 0; \\
\gamma (N^i\zeta , N^j\zeta ) & = & (-1)^i \gamma (w, N^{i+j}w) \, = \, 0.
\end{eqnarray*}
Thus the Gram matrix of $\gamma $ with respect to the basis 
\begin{displaymath}
\{ N^{h+1}\eta = v^0;\, \eta , \, N\eta , \, \ldots , N^h\eta ; 
\, N^h\zeta, - N^{h-1}\zeta , \, \ldots , (-1)^h\zeta  \} 
\end{displaymath}
is given by (\ref{eq-thm5pt3gram}). The matrix of $N$ with respect to the 
above basis in given by (\ref{eq-thm5pt3nil}). \medskip 

Suppose that the nilpotent special type $\underline{\underline{\Delta }}$ is 
decomposable, that is $\underline{\underline{\Delta }} = 
{\underline{\underline{\Delta }}}' + \Delta $, where ${\underline{\underline{\Delta }}}'$ 
is a nilpotent special type and $\Delta $ is a type. Using lemma 2 we may 
assume that ${\underline{\underline{\Delta }}}'$ is indecomposable and has 
the same special height $h+1$ as $\underline{\underline{\Delta }}$. The type 
$\overline{\underline{\underline{\Delta }}}$ induced by the special type 
$\underline{\underline{\Delta }}$ is in case \textbf{a} equal to 
${\Delta }^{\varepsilon}_{h+1}(0)$ and 
in case \textbf{b} is either ${\Delta }_{h+1}(0,0)$ if $h$ is odd or 
${\Delta }^{+}_{h+1}(0) + {\Delta }^{-}_{h+1}(0)$ by construction. 
Because the special height of the special types ${\underline{\underline{\Delta }}}'$ 
and $\underline{\underline{\Delta }}$ are equal to $h+1$, 
we can repeat the construction given in case \textbf{a} or \textbf{b} on the special type 
$\overline{\underline{\underline{\Delta }}'}$ 
induced from the special type ${\underline{\underline{\Delta }}}'$. Hence 
$\overline{\underline{\underline{\Delta }}'}$ in case \textbf{a} is 
${\Delta }^{\varepsilon}_{h+1}(0)$ and in case \textbf{b} is ${\Delta }_{h+1}(0,0)$ 
if $h$ is even or ${\Delta }^{+}_{h+1}(0) + {\Delta }^{-}_{h+1}(0)$ if $h$ is odd.  
Since the nilpotent special types $\underline{\underline{\Delta }}$ and 
${\underline{\underline{\Delta }}}'$ are each uniquely determined by their respective 
induced types, which we have shown are equal, it follows that 
$\underline{\underline{\Delta }} =  {\underline{\underline{\Delta }}}'$. Hence the nilpotent special type $\underline{\underline{\Delta }}$ is indecomposable. \medskip 

This completes the proof of the second and third items of theorem 5 and thereby the 
classification of indecomposable nilpotent special 
types. \hfill $\square $ \medskip %

\noindent \textbf{Proposition 6.} Let $\underline{\underline{\Delta }}$ be 
a special type of special height $h+1$. Suppose that 
\begin{equation}
\underline{\underline{\Delta }} = {\underline{\underline{\Delta }}}' + 
\Delta, 
\label{eq-sec4onenew}
\end{equation}
where ${\underline{\underline{\Delta }}}'$ is an indecomposable 
nilpotent special type of special height $h+1$ and $\Delta $ is 
a type of height less than or equal to $h$. Then the 
decomposition (\ref{eq-sec4onenew}) is unique up to reordering of the 
summands in the type $\Delta $. \medskip 

\noindent \textbf{Proof.} The decomposition (\ref{eq-sec4onenew}) 
yields the decomposition $\overline{\underline{\underline{\Delta }}} 
= {\overline{\underline{\underline{\Delta }}}}{\, '} + 
\overline{\Delta}$ for the induced types. By results of 
\cite{burgoyne-cushman} the induced type 
$\overline{\Delta }$ is the unique sum of 
indecomposable types ${\Delta }_i$ for $1 \le i \le s$ of height $h_i$ 
where $h_1 \le h_2 \le \cdots \le h_{s}$. Therefore  
\begin{displaymath}
\overline{\Delta } 
= {\Delta }_1 + \cdots + {\Delta }_s + {\overline{\Delta }}' .  
\end{displaymath}
Note that $\Delta = \overline{\Delta }$. This follows because any type has 
a representative pair $(N|\widetilde{V}, \widetilde{V})$, where 
$v^0 \not \in \widetilde{V}$. Therefore the projection map 
$\pi $ in the definition of induced type, when restricted to 
$\widetilde{V}$, is a bijective linear map, which defines an 
equivalence of pairs. So $\Delta = {\Delta }_1 + \cdots + 
{\Delta }_s$ is the unique decomposition of $\Delta $ into 
indecomposable types. Using the classification 
of indecomposable nilpotent special types, we see the indecomposable nilpotent 
special type ${\underline{\underline{\Delta }}}'$ in the decomposition 
(\ref{eq-sec4onenew}) is uniquely determined by its induced type 
${\overline{\underline{\underline{\Delta}}}}{\, '}$. Thus the decomposition 
(\ref{eq-sec4onenew}) is unique up to reordering of the summands in 
the type $\Delta $. \hfill $\square $ \medskip 

Combining the results of lemma 2 and 6 we obtain \medskip 
    
\noindent \textbf{Theorem 7.} Every special type may be written as 
a sum of an indecomposable nilpotent special type and a sum of 
indecomposable types. This decomposition is unique up to 
reordering of the indecomposable type summands. \medskip 

Recall that indecomposable types have been classified in 
\cite{burgoyne-cushman}. \medskip 

Theorem 7 solves the problem of finding a unique representative 
for each conjugacy class of elements of ${\oo (V,\gamma )}_{v^0}$ under 
the adjoint action of ${\OO (V, \gamma )}_{v^0}$.

\end{document}